\documentclass{amsart}
\usepackage{amssymb}
\newtheorem{thm}{Theorem}[section]
\newtheorem{cor}[thm]{Corollary}
\newtheorem{lem}[thm]{Lemma}

\theoremstyle{definition}

\theoremstyle{remark}

\numberwithin{equation}{section}

\newcommand{\R}{\mathbb R}

\newcommand{\la}{\lambda}

\begin{document}

\title[Ricci-flat metrics]{A gap theorem for Ricci-flat 4-manifolds}
\author{Atreyee Bhattacharya}
\author{Harish Seshadri}

\address{department of mathematics,
Indian Institute of Science, Bangalore 560012, India}
\email{atreyee@math.iisc.ernet.in}
\email{harish@math.iisc.ernet.in}




\vspace{3mm}

\begin{abstract}
Let $(M,g)$ be a compact Ricci-flat $4$-manifold. For $p \in M$ let $K_{max}(p)$ (respectively $K_{min}(p)$) denote the maximum (respectively the
minimum) of sectional curvatures at $p$. We prove that if
$$K_{max} (p) \ \le \ -c K_{min}(p)$$
for all $p \in M$, for some constant $c$ with $0 \leq c < \frac{2+\sqrt 6}{4}$, then $(M,g)$ is flat.

We prove a similar result for compact Ricci-flat K\"ahler surfaces. Let $(M,g)$ be such a surface and for $p \in M$ let $H_{max}(p)$ (respectively $H_{min}(p)$) denote the maximum (respectively the
minimum) of holomorphic sectional curvatures at $p$. If
$$H_{max} (p) \ \le \ -c H_{min}(p)$$
for all $p \in M$, for some constant $c$ with $0 \leq c < \frac {1+\sqrt 3 }{2}$, then $(M,g)$ is flat.

\end{abstract}


\maketitle

\section{Introduction}
Let $(M,g)$ be a compact Einstein $n$-manifold, $n \ge 4$. One is interested in understanding curvature
conditions which force $g$ to be locally symmetric. More precisely, one seeks pointwise restrictions on curvature (not involving global quantities such as volume and diameter) which imply local symmetry of $(M,g)$. One of the early results in this direction
is due to Tachibana \cite{Tachi}: If $(M,g)$ has nonnegative curvature operator then it is a compact symmetric space. S. Brendle \cite{Bren} generalized Tachibana's result by obtaining the same conclusion under the assumption that $(M,g)$ has nonnegative isotropic curvature. In particular, the result is also true if $(M,g)$ has positive quarter-pinched sectional curvature.

The case of four-dimensional Einstein manifolds has attracted special attention, as it is the lowest dimension
for which the Einstein condition is strictly weaker than that of constant sectional curvature. C. LeBrun and M. Gursky \cite{Gur}
proved that if $(M,g)$ is a compact oriented Einstein $4$-manifold with positive intersection form and nonnegative sectional curvature then $(M,g)$ is isometric to $({\mathbb C}P^2,g_0)$, where $g_0$ is a constant multiple of the Fubini-Study metric. Yang \cite{Yang} proved the following: Let $(M,g)$ be a compact Einstein $4$-manifold normalized so that
$Ric=g$. If the sectional curvatures $K$ satisfy $K \ge 0.1$ then $(M,g)$ is isometric to a round $4$-sphere or to
$({\mathbb C}P^2,g_0)$, $g_0$ as above.

   The above results concern Einstein manifolds with {\it positive} scalar curvature. There are relatively few results when the scalar curvature is zero or negative. The only such result we are aware of is that of Y. T. Siu and P. Yang \cite{Siu}:
   Suppose that $(M,g)$ is a compact K\"ahler-Einstein surface with nonpositive holomorphic bisectional curvature. Let $H_{av}(p)$(respectively $H_{max}(p), H_{min}(p)$) denote the average (respectively maximum, minimum) holomorphic sectional curvature at $p$.  Suppose for some $c < \frac{2}{3(1+\sqrt{\frac{6}{11}})}$, one has
$$H_{av}(p)-H_{min}(p) \le c(H_{max}(p)-H_{min}(p))$$ at every $p \in M$. Then $(M,g)$ is a complex-hyperbolic
space form.

The first result we prove in this paper arose out of an attempt to prove an analogue of the Siu-Yang result for Riemannian $4$-manifolds when the scalar curvature is zero or negative. While the method we employ does not seem to work in the negative case, we have the following result for Ricci-flat $4$-manifolds.

\begin{thm}\label{lk}
Let $(M,g)$ be a compact Ricci-flat $4$-manifold. For $p \in M$ let $K_{max}(p)$ (respectively $K_{min}(p)$) denote the maximum (respectively the
minimum) of sectional curvatures at $p$. If
$$K_{max} (p) \ \le \ -c K_{min}(p)$$
for all $p \in M$, for some constant $c$ with $0 \leq c < \frac{2+\sqrt 6}{4}$, then $(M,g)$ is flat.
\end{thm}

This condition appears natural in light of the following observation (which follows from Lemma \ref{gg} below): If $(M,g)$ is an arbitrary Ricci-flat $4$-manifold, then
\begin{equation}\label{pop}
-\frac{1}{2}K_{min}(p) \le K_{max}(p) \le -2K_{min}(p)
\end{equation}
for any $p \in M$.  Now let $a, c$ be two constants  such that $$ -a K_{min}(p) \le K_{max}(p) \le -c K_{min}(p)$$
for each $p \in M$. Then it follows from the pair of inequalities \eqref{pop} that if $0 \le  c < \frac{1}{2}$ or $a > 2$, then $(M,g)$ must be flat. Note that in Theorem \ref{lk} one assumes that $0 \leq c < \frac{2+\sqrt 6}{4} \sim 1.112$ whereas the same conclusion follows.\\

Our second result is an analogue of the Siu-Yang result for Ricci-flat K\"ahler surfaces:

\begin{thm}\label{lk2}
Let $(M,g)$ be a compact Ricci-flat K\"ahler surface. For $p \in M$ let $H_{max}(p)$ (respectively $H_{min}(p)$) denote the maximum (respectively the
minimum) of holomorphic sectional curvatures at $p$. If
$$H_{max} (p) \ \le \ -c H_{min}(p)$$
for all $p \in M$, for some constant $c$ with $0 \leq c < \frac {1+\sqrt 3 }{2}$, then $(M,g)$ is flat.
\end{thm}

The proof of Theorem \ref{lk} is based on a Bochner-type formula for the curvature tensor of any Einstein $n$-manifold. This formula, which  follows easily from the evolution equation of the curvature tensor under the Ricci flow, actually requires the curvature tensor to be Einstein rather than just harmonic. Given the Bochner formula, one evaluates it at a point where the minimum sectional curvature is attained.  The Laplacian of the sectional curvature function should be nonnegative at this point. However, if the metric is not flat and the pinching condition in Theorem \ref{lk} holds, then the Bochner formula will imply that the Laplacian is negative. To see this, we use a special orthonormal frame for Einstein $4$-manifolds constructed by M. Berger ~\cite{Ber}.  The proof of Theorem \ref{lk2} is similar except that we use the Bochner formula and special orthonormal frame from ~\cite{Siu}.

\section{Proof}
\subsection{Preliminaries:}

The first lemma is an easy consequence of the evolution equation for the curvature tensor under the Ricci flow :
\begin{lem}
Let $(M,g)$ be an Einstein $n$-manifold with $Ric=\lambda g$. We then have
$$ \Delta R + Q(R) =2\lambda R ,$$
where, in an orthonormal frame,
$$Q(R)_{ijkl}=2(B_{ijkl}-B_{ijlk} +B_{ikjl}-B_{iljk})$$
and
$$B_{ijkl}=R_{ipjq}R_{kplq}$$
where the sign convention for $R_{ijkl}$ is such that $R_{ijij}$ is  sectional curvature.
\end{lem}

The next lemma is elementary and well-known:
\begin{lem}
Let $(M,g)$ be a Riemannian $4$-manifold. $(M,g)$ is Einstein if and only if for every $p \in M$ and
$2$-dimensional subspace $V \subset T_pM$, we have
$$K(V)=K(V^\perp)$$
where $K(V)$ denotes sectional curvature of $V$.
\end{lem}

Finally, we have a lemma of Berger \cite{Ber}:
\begin{lem}\label{gg}[M. Berger]
Let $(M,g)$ be an Einstein $4$-manifold and $p \in M$. Then there exists an orthonormal basis $\{e_0,e_1,e_2,e_3\}$ of
$T_pM$ such that the following hold: \\
(1) $K_{01}=  K_{max}(p) $ \\
(2) $K_{03} = K_{min}(p)$ \\
(3) $R_{ijik}=0$ for $j \neq k$ and $0 \le i,j,k \le 3$. \\
(4) $\vert R_{0123} -R_{0231} \vert \ \le \ K_{01}-K_{02}$ \ and  \ $\vert R_{0231} -R_{0312} \vert \ \le \ K_{02}-K_{03}$. \\
\end{lem}
An immediate corollary is 
\begin{cor}\label{one}
Let $(M,g)$ be a Ricci-flat manifold. We have
$$-\frac{1}{2}K_{min}(p) \le K_{max}(p) \le -2K_{min}(p)$$
at every $p \in M$.
\end{cor}

\subsection{$Q(R)$ in dimension 4:}

Let $(M,g)$ be an Einstein $4$-manifold with $Ric=\lambda g$. In this section we find an expression for the ``sectional curvatures" of the $(0,4)$-tensor $Q(R)$:
\begin{align}
Q(R)_{ijij} &= 2(B_{ijij} - B_{ijji}+ B_{iijj} - B_{ijji}) \\ \notag
            &= 2(R^2_{ipjq}+R_{ipiq}R_{jpjq} - 2R_{ipjq}R_{jpiq})
\end{align}
Let $\{e_0,e_1,e_2,e_3 \}$ be an orthonormal basis as in Lemma \ref{gg}: We then have
\begin{align} \notag
\frac{1}{2}Q(R)_{0303} &= K_{03}^2 +R_{0132}^2+R^2_{0231} -2R_{0132}R_{3102}-2R_{0231}R_{3201} +K_{01}K_{31}+K_{02}K_{32} \\ \notag
            &= K_{03}^2+ 2K_{01}(\lambda -K_{01}-K_{03}) + R_{0132}^2+R^2_{0231}-4R_{0132}R_{3102}.
\end{align}
In the second equality we have used the following chain of equalities:
Since
$$K_{31}=K_{02} \ \ {\rm and} \ \ K_{32}=K_{01}$$
we have
\begin{equation}
 K_{01}K_{31}+K_{02}K_{32} = 2K_{01}K_{02}= 2K_{01}(\lambda -K_{01}-K_{03}). \\ \notag
\end{equation}

Therefore
$$\frac{1}{2}Q(R)_{0303} = K_{03}^2+2(\la-K_{03})K_{01}-2K_{01}^2 + q$$
where
$$q=R_{0132}^2+R^2_{0231}-4R_{0132}R_{3102}.$$

\vspace{3mm}

\subsection{Estimates for $q$:} \ \ \
\vspace{2mm}

In this section assume that $(M,g)$ is a Ricci-flat $4$-manifold which is not flat. For $p \in M$ assume that the basis of
$T_pM$ coming from Lemma \ref{gg} satisfies $K_{01}=K_{max}=1, \ K_{03}=K_{min}= - \delta$. One then has
$K_{02}=\delta -1$. Also it follows from Corollary \ref{one} that $$\frac{1}{2} \le \delta  \le 2.$$

First note that
\begin{align}\notag
 q &= R_{0132}^2+R^2_{0231}-4R_{0132}R_{3102} \\ \notag
   &=R_{0123}^2+R_{0231}^2 +4 R_{0123}R_{0231} \notag
\end{align}
\vspace{2mm}

Let $x=R_{0123}$, $y=R_{0312}$ and $z=R_{0231}$. Then
\begin{equation}\label{po}
\vert x-z \vert \le 2 -\delta, \ \ {\rm and} \ \ \vert z-y \vert \le 2 \delta -1.
\end{equation}
Since
$$x+y+z=0$$
(\ref{po}) gives
\begin{equation}\label{po1}
\vert x-z \vert \le 2 -\delta, \ \ {\rm and} \ \ \vert x+2z \vert \le 2 \delta -1.
\end{equation}
\vspace{2mm}

We want to find the global maximum and minimum of
$$q(x,z)= x^2+z^2+4xz$$
on the polygonal region $D$ defined by the inequalities (\ref{po1}) in $\R ^2$ (with coordinates $(x,z)$).
\vspace{2mm}

The four line segments bounding $D$ are given by
$$ z=x \pm (2-\delta), \ \ z = \frac{-x \pm (2\delta-1)}{2}$$
and the vertices of $D$ are
$$p_1=(1, \ \delta-1), \ \ p_2=-p_1, \ \ p_3= \Bigl (\frac {4\delta -5}{3}, \ \frac { \delta+1}{3} \Bigr ), \ \ p_4=-p_3.$$
\vspace{2mm}

The values of $q$ at these corner points are

\begin{equation}\label{an}
 q(p_1)=q(p_2)=  \delta^2 +2 \delta-2, \ \ \ q(p_3)=q(p_4)= \frac {33\delta^2 - 42 \delta +6}{9} = \frac{11\delta^2 - 14\delta + 2}{3}
\end{equation}

Note that $q$ cannot have a local maximum or minimum in the interior of $D$ as
$$ {\rm Hessian}(q)=  \left( \begin{array}{cc}
2 & 4  \\
4 & 2  \end{array} \right)  $$
is neither positive nor negative definite.
If $q$ has a critical point $(x,z)$ in the interior of any line segment comprising the boundary of $D$, then
$$ \langle \nabla q (x,z), X \rangle =0$$
where $X=(1,m)$ is a nonzero tangent vector to the line at $(x,z)$.
This gives
$$ x+2z +m(2x+z)=0.$$
We then have two possibilities:
\begin{equation}\label{mk}
m=1, \ \ \  x=-z= \pm \frac{2 -\delta }{2} \ \ \ {\rm and} \ \ \ q= -\frac{(2 -\delta)^2}{2}
\end{equation}
or
\begin{equation}\label{mk2}
m=-\frac{1}{2} \ \ \  z=0, \ \ x = \pm (2\delta-1) \ \ \ {\rm and} \ \ \ q=(2\delta-1)^2.
\end{equation}

Combining (\ref{an}), (\ref{mk}) and (\ref{mk2}), we see that
\begin{equation}\notag
q \ \ge \ min \ \Bigl \{\delta^2 +2 \delta-2, \ \ \frac {11\delta^2 - 14 \delta +2}{3}, \ \ -\frac{(2 - \delta)^2}{2} \Bigr \}
\end{equation}
We note the following facts: First,  \ $min \ \Bigl \{ (\delta^2 + 2\delta- 2), \  \frac {(11\delta^2 - 14 \delta +2)}{3} \Bigr \}  \ \ge \ -\frac{(2 -\delta)^2}{2}$ for all values of $\delta$. Second, the critical points $(\pm \frac{2 - \delta}{2}, \mp \frac{2 - \delta}{2})$ corresponding to the critical value $-\frac{(2 - \delta)^2}{2}$ of $q$ belong to the polygonal region $D$ if and only if $\delta \ \ge \ \frac{4}{5} = 0.8$. Hence

\begin{equation}\label{pk}
q \ \ge \ \ \ \begin{cases} -\frac{(2-\delta)^2}{2} \ \ \ \ \ \ \ \ \ \ \ \ \ \ \ \ \ \ \ \ \ \ \ \ \ \ \ \ \ \ \ \ \ \  {\rm if} \ \ \delta \ \ \ge \  \frac{4}{5};\\
min \ \Bigl \{\delta^2 +2 \delta-2, \ \ \frac {11\delta^2 - 14 \delta +2}{3}\Bigr \} \ \ {\rm otherwise.}
\end{cases}
\end{equation}

\vspace{3mm}

\subsection{Sign of $Q(R)$:} \ \ \
\vspace{2mm}

When $(M,g)$ is a Ricci-flat $4$-manifold and $\{e_0,e_1,e_2,e_3 \}$ as in Lemma \ref{gg}, one has
\begin{align}\notag
Q(R)_{0303} &= K_{03}^2+2(\la-K_{03})K_{01}-2K_{01}^2 + q \\ \notag
            &= \delta^2+2\delta -2 +q \notag
\end{align}
We now use the lower bound (\ref{pk}): \vspace{2mm}

If \ $q  \ \ge \  -\frac {(2 - \delta)^2}{2}$, then
$$ Q(R)_{0303}  \ \ge \ \frac{1}{2} \Bigl (\delta^2 + 8\delta - 8 \Bigr)$$
which implies that
\begin{equation}
Q(R)_{0303} > 0 \ \ \ {\rm if} \ \ \  \delta >  2(\sqrt{6} - 2) \sim 0.8989
\end{equation}

\vspace{3mm}

In conclusion, we have
$$ Q(R)_{0303} > 0 \ \ \ {\rm if} \ \ \ \delta^{-1} <  \frac{1}{2(\sqrt{6} - 2)} = \frac{(\sqrt{6} + 2)}{4} \sim 1.112 .$$
Although if $\delta < \frac{4}{5},$ then nothing can be concluded about the sign of $Q(R)_{0303}$ from the above calculations.

\subsection{Proof of Theorem \ref{lk} } \

\vspace{2mm}

Suppose $(M,g)$ is a compact Ricci-flat 4-manifold such that for all $p \in M$, $K_{max}(p) \le -c K_{min}(p)$ where $0 < c <  \frac{(\sqrt{6} + 2)}{4}$.
We will assume that $(M,g)$ is not flat and get a contradiction.
Let $p$ be a point in $M$ where the minimum sectional curvature (over all  2-dimensional subspaces of tangent spaces) $K_{min}(p) \neq 0$ is attained. Rescale $g$ so that $K_{max}(p)=1$. If $\delta: = -K_{min}(p)$ then $\delta^{-1} <  \frac{(\sqrt{6} + 2)}{4}$ by hypothesis.

  Let $\{e_0,e_1,e_2,e_3\}$ be a basis of $T_pM$
as in Berger's Lemma. Choose $r$ less than the injectivity radius at $p$ and extend the basis to an orthonormal frame
$\{E_0,E_1,E_2,E_3\}$ on $B_p(r)$ by parallel translating $e_0,e_1,e_2,e_3$ along radial geodesics (starting at $p$). One then has
$$ \nabla _{E_i}E_j (p)= \nabla_{E_i}\nabla _{E_i}E_j(p) =0$$
for all $i,j$ and hence
$$  \Delta (R(E_0,E_3,E_0,E_3))(p)=(\Delta R) (e_0,e_3,e_0,e_3) $$
Since $p$ is a minimum for the function $x \mapsto R(E_0(x),E_3(x),E_0(x),E_3(x))$, we have
$$\Delta (R(E_0,E_3,E_0,E_3))(p) \ \ge \ 0.$$
On the other hand, $(\Delta R) (e_0,e_3,e_0,e_3) = -Q(R)(e_0,e_3,e_0,e_3) < 0 \ \  {\rm if} \ \ \delta^{-1} <  \frac{(\sqrt{6} + 2)}{4}$. This contradiction completes the proof. \hfill $\square$
\vspace{3mm}

{\bf Remark:} Instead of considering the point of minimum for sectional curvatures we can also consider the point
of maximum. In this case we work with the equation
$$(\Delta R)_{0101} = -Q(R)_{0101}.$$
Now
$$\frac{1}{2}Q(R)_{0101}=  K_{01}^2-2K_{01}K_{03}-2K_{03}^2 + q$$
with
$$q = R_{0231}^2+R^2_{0312}+4R_{0231}R_{0312}.$$
In order to estimate $q$ one works with the region defined by the inequalities
$$ \vert y+2z \vert \le 2 -\delta, \ \ {\rm and} \ \ \vert z-y \vert \le 2 \delta -1$$
and sees that
$$q \ \le \ max \ \Bigl \{  \frac {2 \delta^2 - 14 \delta +11}{3}, \ \ (2-\delta)^2, \ \ (1 + 2\delta - 2\delta^2) \Bigr \} = (2-\delta)^2 \ \ {\rm if \ \ \delta \ge 1}.$$
Proceeding as before, one finds that
$$Q(R)_{0101} < 0 \ \ \ {\rm if} \ \ \  \delta > \sqrt{6}-1 \sim 1.45 $$ or equivalently,
$$Q(R)_{0101} < 0 \ \ \ {\rm if} \ \ \  \delta^{-1} < \frac{\sqrt{6} + 1}{5} \sim 0.69. $$
Hence one gets a weaker restriction on $\delta$ (respectively $\delta^{-1}$) by considering the minimum. \vspace{3mm}

\subsection{K\"ahler-Einstein surfaces} \ \
\vspace{2mm}
This section is based on the computations of Siu-Yang ~\cite{Siu} and we refer the reader to it for details. \vspace{2mm}

Let $(M,J,g)$ be a compact K\"ahler-Einstein surface, with $J$ denoting the almost-complex structure, and let $p \in M$. Then one has
\begin{equation}\label{ke}
\frac{1}{3}(H_{max}(p)- H_{min}(p)) \ \le \ (H_{av}(p) - H_{min}(p)) \ \le \ \frac{2}{3}(H_{max} - H_{min}(p))
\end{equation}
where $H$ denotes hololomorphic sectional curvature. If we assume that assume that $(M,J, g)$ is Ricci-flat, then  
$H_{av}=0$ and the pair of inequalities \eqref{ke} reduce to
\begin{equation}\label{ke2}
-\frac{1}{2}H_{min}(p) \ \le H_{max}(p) \ \le -2H_{min}(p).
\end{equation}
Let $a, c$ be constants such that  $$ -a H_{min}(p) \le H_{max}(p) \le -c H_{min}(p)$$
for each $p \in M$. Then as in the previous case of Ricci-flat Riemannian $4$-manifolds, it follows from the inequalities \eqref{ke2} that if $0 \le  c < \frac{1}{2}$ or $a > 2$, then $(M,g)$ is flat.\\

Now let $\{e_1, \ e_2\}$ be an unitary frame of $T^0_pM$, the holomorphic tangent space at $p$. Note that
$$e_1= \frac {1}{\sqrt 2} (u -\sqrt{-1}J(u)), \ \ \ e_2= \frac {1}{\sqrt 2} (v -\sqrt{-1}J(v))$$
for an orthonormal set $\{u,J(u),v,J(v)\} \subset  T_pM$.
Also
$$ R_{1 \bar 1 1 \bar 1} = K(u,J(u))=H(e_1), \  \ \ R_{1 \bar 1 2 \bar 2} = K(u,v)+K(u,J(v))$$
where $H(e_1)$ denotes the holomorphic sectional curvature of the complex line spanned by $e_1$.

Suppose
$$ Ric = \lambda g.$$ Let $R_{1 \bar 1 1 \bar 1} = H(e_1)$ be critical among all holomorphic sectional curvatures at $p$.
Then $R_{\alpha \bar \beta \gamma \bar \delta}$ vanishes at $p$ whenever precisely three of $\alpha, \beta, \gamma, \delta$ are equal and we have the following Lemma due to Siu and Yang \cite{Siu}
\begin{lem}\label{two}{\rm (Proposition 2, ~\cite{Siu})}
$$(\Delta R)_{1 \bar 1 1 \bar 1} = -AR_{1 \bar 1 2 \bar 2} + \vert B \vert^2,$$
where
$$ A= 2R_{1 \bar 1 2 \bar 2} - R_{1 \bar 1 1 \bar 1}, \ \ \ B = R_{1 \bar 2 1 \bar 2}.$$
\end{lem}

Also
$$H_{av}= \frac {2}{3} \lambda = \frac {2}{3}(R_{1 \bar 1 1 \bar 1} + R_{1 \bar 1 2 \bar 2}).$$

Assume further that $R_{1 \bar 1 1 \bar 1} = H_{min}(p)$.
Then $$A = -3H_{min}(p) \ \ {\rm and} \ \ \vert B \vert = 2 H_{max}+H_{min}.$$
Using Lemma \ref{two} we obtain
\begin{equation}
(\Delta R)_{1 \bar 1 1 \bar 1}  =  -3 H_{min}^2 + ( 2 H_{max}+H_{min})^2.
\end{equation}
Therefore
$$(\Delta R)_{1 \bar 1 1 \bar 1} <0 \ \ \ {\rm if} \ \ \     H_{max} \le -c H_{min}$$
for any $c$ satisfying
$$ 0 \le  c <  \frac {\sqrt 3 +1}{2} \sim 1.366.$$
\vspace{3mm}

{\bf Remark:} Instead of considering the point of minimum for holomorphic sectional curvatures consider the point
of maximum. Assume that $R_{1 \bar 1 1 \bar 1} = H_{max}(p).$ Then we obtain the following expression
$$(\Delta R)_{1 \bar 1 1 \bar 1} = -AR_{1 \bar 1 2 \bar 2} + \vert B \vert^2,$$
where
$$ A= - 3H_{max}(p), \ \ \ \vert B \vert = -(H_{max}(p) + 2H_{min}(p)).$$
Thus in this case
$$(\Delta R)_{1 \bar 1 1 \bar 1}  =  -3 H_{max}^2 + ( 2 H_{min}+H_{max})^2 > 0 \ \ {\rm if} H_{max} \le -c H_{min}$$
for any $c$ satisfying
$$ 0 \le  c <  \sqrt 3 - 1 \sim 0.73.$$
which is a stronger restriction on $c$ than the one obtained by considering the minimum.

\end{document}